\documentclass{article}
\usepackage{graphicx}
\usepackage{amssymb, amsfonts, amsbsy, latexsym, amsmath}
\usepackage{tikz}
\usepackage{color}
\usepackage{colortbl}
\usepackage{url}
\usepackage{hyperref}
\usepackage{bezier}

\begin{document}


\begin{titlepage}
\title{A Page-Rank-like Approach to Optimal Placement of Charging Stations in a Warehouse} 
\author{Hans-Georg Stark, Mustafa Jelibaghu, Katrin Tschirpke, Michael Eley\\
Technische Hochschule Aschaffenburg, W\"urzburger Straße 45, 63743 Aschaffenburg\\
E-Mail: {\tt FirstName.LastName@th-ab.de}\\ 
\date{\today}} 
\end{titlepage}
\maketitle

\begin{abstract}
\noindent In this paper we describe an approach to the problem of finding optimal positions for charging stations (CS) in a warehouse, where a fleet of electrical industrial trucks/forklifts is employed. \\
Our procedure is motivated by Google's Page-Ranking. The graph model underlying our method is easily extensible from simple distance based criteria, relevant for choosing optimal CS-positions, to more complex criteria taking into account, e.g., the ''state of charge'' (SOC) of the individual truck's battery. 
\end{abstract}

MSC-Codes: 05C90, 94C15 

\section{Introduction}
We consider a fleet of electrically driven industrial trucks/forklifts on a warehouse site. The road map of this site is known and has a graph structure consisting of vertices and edges (see below). The task is, to select a certain number of those vertices as locations of charging stations (CS-locations), based on historically recorded traffic data and/or simulated data, provided such a simulation environment is given. 

This problem and variants thereof have been treated mostly in the context of urban electrical mobility, see, e.g., \cite{EVCSP, EVCSP2, EVCSP3}.

In \cite{EVCSP}, an optimization problem is formulated as a combinatorial problem, which is solved by integer linear programming techniques. The algorithm determines a minimum number of strategically chosen CS-locations, given - among others - the range of electrical vehicles.
In \cite{EVCSP2} an approach is described, which is based on statistical properties of spatio-temporal vehicle distributions, which are evaluated by machine learning methods. In \cite{EVCSP3}, the Electric Vehicle Charging Station Placement Problem (EVCSPP) is formulated based on some a priori-data like local charging station demand requirements and average range of electrical vehicles. The problem again is formulated as a constrained optimization task and several solution methods are evaluated. 

The approach described in this paper is in a sense simpler, since a warehouse roadmap is less complex than urban or national road networks, moreover we have access to truck position data reflecting typical transport situations, which might be enhanced by warehouse simulations (see below). Moreover, we formulate the relevance of a road-map-vertex as a potential CS-location by borrowing and adapting Google's PageRanking \cite{PRPatent} for our purpose. In the next section we describe our approach, in Sect. \ref{sec:Res} we discuss results and sketch future developments.

\section{Graphs, Page-Ranking and Optimization \label{sec:GPO}}
In this section we describe our approach to formulate and solve the CS-position-problem in a warehouse. It follows a 3-step-scheme.
\begin{enumerate}
    \item \textit{Data acquisition}. We assume that ''typical'' motion data of a fleet of industrial trucks or forklifts are available. They may have been recorded on the warehouse in some representative time period (this is the case in our evaluation). On the other hand they could have been generated via simulation, provided a ''digital twin'' of warehouse and trucks does exist.
    \item \textit{Graph model}. Potential CS-locations are vertices of a warehouse-roadmap. They are linked in a suitable manner to the recorded truck data in order to quantify the suitability of those vertices as CS-locations.\\
    This procedure requires a special graph model explained in the next section. The ''quantification of suitability'' is described in Sect. \ref{sec:PageRank}.
    \item \textit{Optimization}. After each vertex has been labelled with its suitability rating, the most important ones have to be chosen. Since it makes no sense to choose CS-locations which are in immediate neighborhood to each other, proper constraints have to be formulated in order to ensure a meaningful distribution of the calculated CS-locations on the warehouse. This leads to an optimization problem, which is described in Sect. \ref{sec:ILP}.
\end{enumerate}

\subsection{Graphs \label{sec:Graphs}}
Graph theory has been applied extensively in chemistry, operational research, engineering, social networking, logistics and computing (see \cite{stankovic}, \cite{Shuman1}, \cite{Shuman}, for a recent review refer to \cite{ortega}). 

In mathematical terminology a graph is defined as a triple:

$$
G = \lbrace V, E, W \rbrace
$$
 
\noindent Here $V =\lbrace v_1, v_2, ..., v_N\rbrace$ is a finite set of $N$ ''vertices'', $E \subseteq V \times V $ is a set of edges defined as order pairs $(n, m)$ connecting $v_n$ to $v_m$ and the weight function $W : E \rightarrow \mathbb{R}$ assigns the weight value $w_{nm}$ to edge $(n,m)$ with the understanding that $w_{nm}=0 \Longleftrightarrow (n,m)\notin E$. 

Thus the weight matrix $A^W$ with elements $A^W_{nm}=w_{nm}$ completely describes the connectivity of the graph and undirected graphs are characterized by $w_{nm}=w_{mn}$, i.e., by a symmetric weight matrix $A^W$. When this is not the case, we say that the graph is directed.

Obviously it is useful to represent edge $(n,m)$ by an arrow, linking vertex $v_n$ to vertex $v_m$, with tip at $v_m$ if the graph is directed. For undirected graphs (i.e. if $A^W$ is symmetric) the link usually is represented by an undirected line segment or by a double-tip line segment.

In this paper we shall use a special graph structure. The set of $N$ vertices is split in two subsets as follows: $V=(V_C,V_T)$, where $V_C=\lbrace v_1,\ldots,v_{n_C} \rbrace$ and $V_T=\lbrace v_{n_C+1},\ldots,v_{N} \rbrace$. The elements of $V_C$ form an undirected graph, the vertices from $V_T$ are not linked with each other, but they may be linked to one or more vertices from $V_C$. Thus in this case the weight matrix has the structure
\begin{equation}
A^W=\left(
\begin{array}{cc}
  A^{W_C}   &  0^{CT}\\
  A^{W_{TC}}   & 0^{TT}
\end{array}
\right).
    \label{eq:WeightStructure}
\end{equation}
In Eq. \ref{eq:WeightStructure} the symmetric matrix $A^{W_C}$ is $n_C \times n_C$ and represents the undirected graph connecting the vertices from $V_C$. $A^{W_{TC}}$ is $(N-n_C) \times n_C$ and contains the edges from vertices in $V_T$ to vertices in $V_C$. The matrices $0^{CT}$ and $0^{TT}$ contain zeros and are of dimension $n_C\times (N-n_C)$ and $(N-n_C)\times (N-n_C)$, respectively. They guarantee that vertices from $V_C$ are not linked to vertices from $V_T$ and vertices from $V_T$ are not linked with each other.

\begin{figure}[h]\centering
\includegraphics[width=3cm]{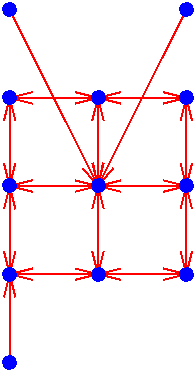}
\caption{An example graph of structure \eqref{eq:WeightStructure}.}%
\label{fig:WeightStructure}%
\end{figure}

An example is shown in Fig. \ref{fig:WeightStructure}. The vertices of $V_C$ form a rectangular $3\times 3$-array, they constitute an undirected graph (here the edges are represented as double-arrows as indicated above). The set $V_T$ consists of the three vertices outside this grid. They are connected to vertices from $V_T$ as indicated.

\begin{figure}[h]\centering
\includegraphics[width=3.8cm]{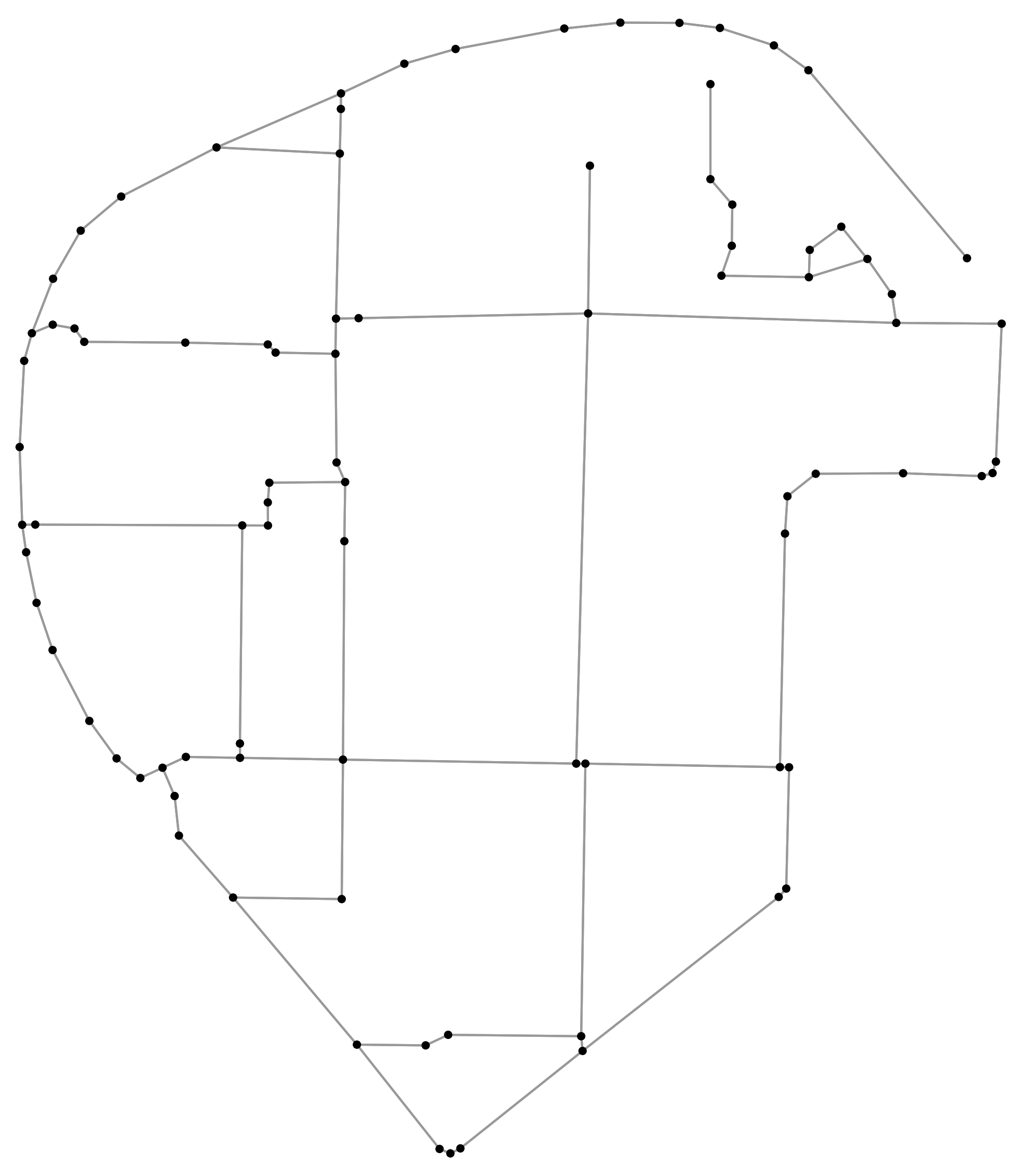}\includegraphics[width=5cm]{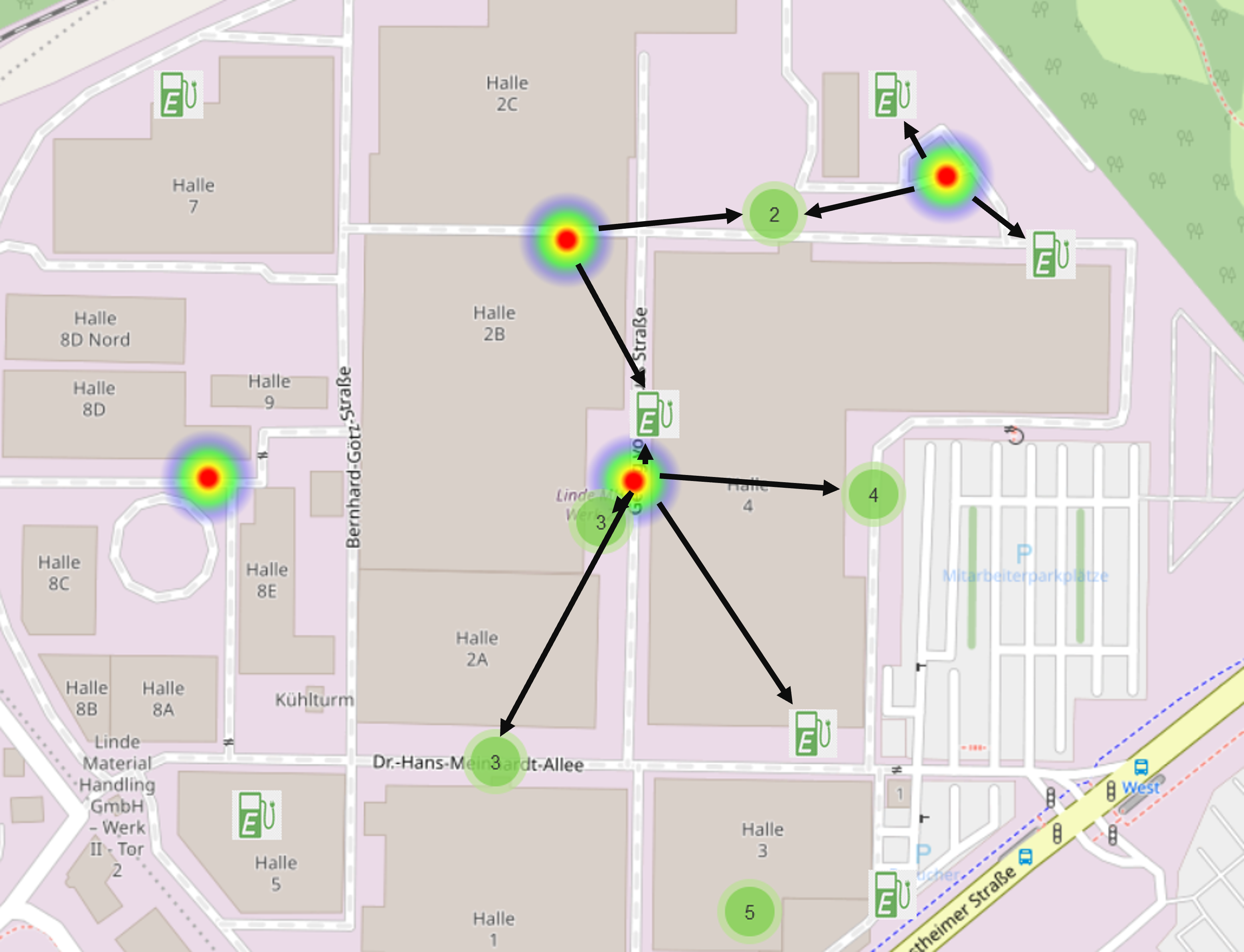}
\caption{Potential CS-Locations $(V_C,A^{W_C})$ (left), partial view of CS-graph with some connected truck vertices $(V,A^W)$.}%
\label{static_dynamic_net}%
\end{figure}

Now we apply this splitting concept to the CS-positioning problem. The sub-graph $(V_C,A^{W_C})$ corresponds to a warehouse-roadmap, whose vertices are potential CS-locations (cf. Fig. \ref{static_dynamic_net}, left). Enhancing this graph with the truck vertices $V_T$, which can be connected to vertices from $V_C$ we obtain the full graph $(V,A^W)$; a partial view is shown in Fig. \ref{static_dynamic_net}, right, the trucks are marked by colored circles.

The rules for connecting trucks to CS-locations are described in the next section.

\subsection{PageRank Algorithm \label{sec:PageRank}} 
Following the procedure sketched above, we turn now to a method for rating the location candidates $v_k \; (k=1,\ldots,n_C)$.

Google's PageRank algorithm \cite{PRPatent} in fact computes the importance of a given web site (vertex) based on the importance of other web sites linked to the given one. In a very simplified fashion the algorithm goes as follows:
\begin{enumerate}
    \item Consider the internet as a directed graph $(V,E,W)$ in the sense defined above. $V=\lbrace v_1,\ldots,v_N \rbrace$ denotes the set web sites, moreover we consider the particular form $W : E \rightarrow \lbrace 0, 1 \rbrace$ (i.e., either there is a link connecting $v_n$ to $v_m$ or not).
    \item Denoting with $p_k$ the importance of $v_k$ the above ''importance idea'' leads to 
    \begin{equation}
    p_k=\sum_{m\ne k} w_{mk} \frac{1}{l_m} p_m,
        \label{eq:PR}
    \end{equation}
    where $l_m=\sum_{n\ne m} w_{mn}$ is the number of outgoing links from $v_m$.
\end{enumerate}
Collecting $p_k$ in a vector $\Vec{p}=(p_1,\ldots,p_N)$, \eqref{eq:PR} can be written as an Eigenvalue-equation:
\begin{equation}
\Vec{p}=M\cdot \Vec{p}
    \label{eq:PR2}
\end{equation}
with the $N\times N$-''Google matrix'' $M$ such that
\begin{equation}
    M_{km}=\frac{w_{mk}}{\sum_{n\ne m} w_{mn}}.
    \label{eq:M}
\end{equation}
This Eigenvalue-equation may be solved by an iterative procedure ${\vec{p}}_{}^{l+1}=M\cdot {\vec{p}}_{}^{l}$ \cite{larypage}.

It turns out that this procedure simplifies considerably, if we apply it to $V=(V_C,V_T)$ with CS-locations $V_C=\lbrace v_1,\ldots,v_{n_C} \rbrace$, truck vertices $V_T=\lbrace v_{n_C+1},\ldots,v_{N-n_C} \rbrace$ and the weight matrix $A^W$ given by \eqref{eq:WeightStructure} as described above. If we assume $p_k=1 \; (k=n_C+1, \ldots, N)$, i.e., all truck vertices have importance $1$, and if we, moreover, shall employ \eqref{eq:PR} only for computing the importance (the rating) $p_k \; k=1,\ldots,n_C$ of the CS-locations and if, finally, we consider in \eqref{eq:PR} only weight contributions from $A^{W_{TC}}$, i.e., links from trucks to CS-locations, the equation reduces to

    \begin{equation}
    p_k=\sum_{m=n_C+1}^N w_{mk} \frac{1}{\sum_{n\ne m} w_{mn}}  \; (k=1,\ldots,n_C)
        \label{eq:PR3}
    \end{equation}
This is a direct non-iterative computation scheme for the CS-location-ratings, which can be efficiently implemented. Note that $\sum_{n\ne m} w_{mn}$ is the $m$-th row sum of $A^W$, which can be computed once, since it does not depend on $k$. To complete the computation scheme, we define now $w_{mk} \; (m>n_C, k\le n_C)$, linking trucks to CS-locations (cf. Fig. \ref{fig:WeightStructure}, right).

Denoting with $d_{mk}$ the Euclidean distance from truck $v_m$ to CS-location $v_k$ the idea is, to have a decreasing edge weight with increasing distance and to remove the edge, if the distance is larger than a certain threshold $t$. Thus a candidate formula reads
	\begin{equation}
    w_{mk} = 
    \begin{cases}
      \frac{1}{1 + d_{mk}}  & \text{if}\ d_{mk} \leq t \\
      0 & \text{else}
    \end{cases}	 \label{eq:wtcs}   
	\end{equation}
As mentioned at the beginning of Sect. \ref{sec:GPO} we need for our evaluation recorded or simulated truck position data from which the above distances and edge weights may be obtained. This implies that a graph like the right part of Fig. \ref{static_dynamic_net} is actually a ''still frame'' belonging to a movie of such graphs. Therefore all quantities above have an additional time label $\lbrace t_l \rbrace$ such that edge weights depend on this label and we obtain via \eqref{eq:PR3} time-dependent CS-position ratings $p^{t_l}_k  \; (k=1,\ldots,n_C)$. Then the final CS-position ratings are obtained by accumulation:
\begin{equation}
p_k=\sum_{t_l} p^{t_l}_k  \; (k=1,\ldots,n_C).
    \label{eq:PR4}
\end{equation}

\subsection{Integer Linear Program \label{sec:ILP}}
After having labelled the vertices of $V_C$, i.e., the CS-position candidates with ratings $p_k$ as described above, the next step is to choose positions with highest rankings. Moreover we want to bound the maximum number of those positions and a ''clustering'' should be avoided by guaranteeing a minimum distance between neighboring selected CS-positions. This leads to a constraint optimization problem, which we formulate in this section. Before doing so, we give some preparatory remarks.

Starting from the CS-Roadmap $(V_C, A^W_C)$ we assume that the nonzero-weight values $A^{W_C}_{ik}$ associated with an edge $(i,k)$ represent the distance of $v_i$ to some connected $v_k$. From the weight matrix the distance matrix $D_{mn} \; m,n=1,\ldots,n_C$ may be derived, where $D_{mn}$ is the length of the shortest path from $v_m$ to $v_n$. This is a standard task and may be solved by Dijkstra's algorithm \cite{SD}. Moreover, we need a ''connectivity matrix'' indicating, if the shortest path connecting two vertices exceeds some threshold $R$:

\begin{equation}
A_{mn} = \begin{cases}
     0  & D_{mn} \le R \\
     1  & D_{mn} > R
   \end{cases} \; (m,n=1,\ldots n_C) \\
  \label{connectivity}
\end{equation}

Next we denote by $X$ the set of binary sequences of length $n_C$:
$$X=\lbrace \Vec{x}=(x_1,\ldots,x_{n_C}) \, | \, x_k \in \lbrace 0,1 \rbrace, \, k=1,\ldots,n_C \rbrace$$
Given some $\Vec{x}\in X$, note, that the quantity $\sum_{k=1}^{n_C}x_kp_k$ sums up only those $p_k$-values, where $x_k=1$. Thus, varying $\Vec{x}$ over $X$ we obtain the accumulated rank values over all possible subsets of $V_C$.

Note furthermore, that for a given $\Vec{x}\in X$ the inequality $x_m+x_n \le 1+A_{mn}$ ensures, that the corresponding selected vertices $v_m$ and $v_n$ have a minimum shortest path distance of $R$ (cf. \eqref{connectivity}). If we, moreover, bound the number of selected CS-locations by $k$, the constraint optimization problem mentioned at the beginning of this section reads
\begin{equation}
    \Vec{x}=\mbox{argmax}_{\Vec{y}\in X} \sum_{k=1}^{n_C}y_kp_k
    \label{eq:Min}
\end{equation}
such that 
\begin{equation}
    x_m+x_n \le 1+A_{mn}, \; m,n=1,\ldots,n_C
    \label{eq:Constr1}
\end{equation}
and
\begin{equation}
    \sum\limits_{m=1}^{n_C} x_m \le k
    \label{eq:Constr2}
\end{equation}

This problem belongs to the class of Integer Linear Programming (ILP) \cite{ILP}, for which standard solution methods are available.

\section{Result, Discussion and Future Work \label{sec:Res}}
\begin{figure}[h]\centering
\includegraphics[width=4cm]{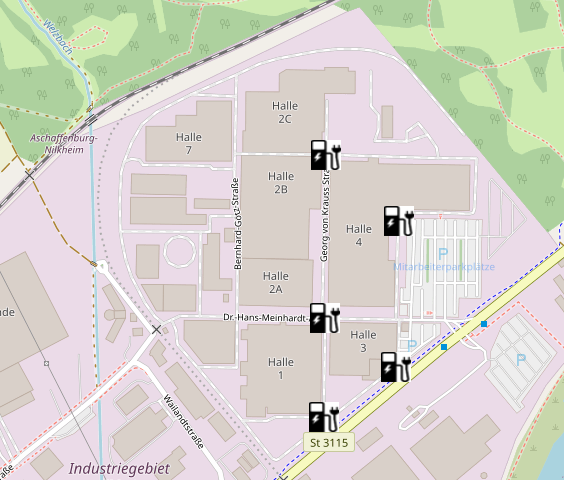}\includegraphics[width=2.7cm]{linde.png}
\includegraphics[width=4cm]{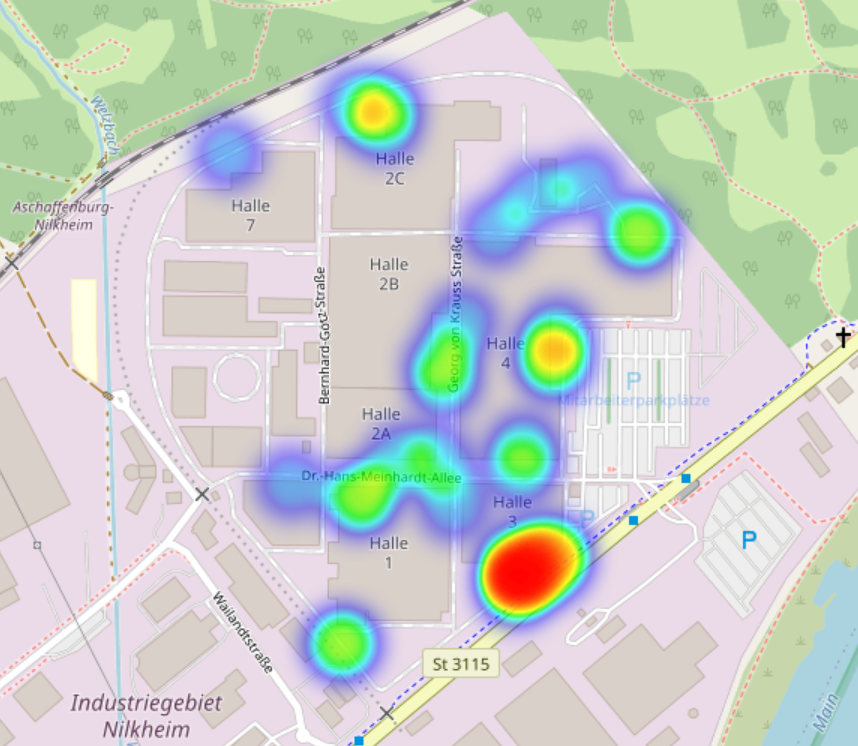}
\caption{Left: Example result ($k=5$), middle: $(V^C,A^{W_C})$, right: Traffic-Heatmap}%
\label{opl_result}%
\end{figure}

The above procedure has been tested at the production site of Linde Material Handling GmbH, a large German forklift manufacturer. The potential CS-locations $(V_C,A^{W_C})$ (Fig. \ref{opl_result}, middle) have been generated via Open Street Map, the truck position data time series have been provided by the company. The Integer Linear Program has been solved with IBM ILOG CPLEX \cite{CPLEX}

The positions selected by the algorithm are shown in Fig. \ref{opl_result}, left. The right part shows for comparison the result of a heat map analysis, showing, which part of the warehouse is heavily utilized, based on the recorded motion data; for a related approach, see also \cite{heatmapanalyse}. It is obvious that the selected CS-positions relate to intense utilization of the factory premises, which makes sense. 

Some improvements of the described procedure are obvious and shall be included in further versions. 

Taking into account that the forklifts are moving on warehouse roadmap it will be more realistic to replace in the  weight function \eqref{eq:wtcs} the Eucidean distance $d_{mk}$ with the shortest path distance $d^{SP}_{mk}$, which again may be computed with Dijkstra's algorithm.\\
Secondly, generalizing the weight structure $A^{W_{TC}}$ from \eqref{eq:WeightStructure} such that the forklift truck's battery state (the ''State of Charge'' (SOC)) is also taken into account would considerably improve the CS-selection procedure, because the algorithm's sensitivity then not only relies on distances from trucks to CS-locations but also on the individual truck's necessity, to reload the battery. This requires either the inclusion of SOC-data in the fleet recording, or - if the data are obtained via simulation - the inclusion of a reliable battery model in the simulation. This is ongoing work.
\section{Acknowledgement}
The work has been supported by the Bavarian Ministry of Economic Affairs, Regional Development and Energy (Project KAnIS, Grant No. DIK-1910-0016//DIK0103/01). The authors are grateful to Linde Material Handling GmbH, the industrial partner in this project, for providing truck data and essential advice. 

\bibliographystyle{plainurl}
\bibliography{references}

\end{document}